\pdfoutput=1
\documentclass[journal]{IEEEtran}
\usepackage{cite}

\usepackage{amsmath}

\usepackage{algpseudocode}

\ifCLASSINFOpdf
  \usepackage[pdftex]{graphicx}
\else
   \usepackage[dvips]{graphicx}
\fi
\ifCLASSOPTIONcompsoc
 \usepackage[caption=false,font=normalsize,labelfont=sf,textfont=sf]{subfig}
\else
\usepackage[caption=false,font=footnotesize]{subfig}
\begin{document}
%
\title{Exact Analysis of Synchronizability for Complex Networks using Regular Graphs}
%
%
%

\author{Sateeshkrishna Dhuli,~\IEEEmembership{Student Member,~IEEE,}
        {Y.N.Singh,~\IEEEmembership{Senior Member,~IEEE}}
}

%
%

\markboth{}%
{Shell \MakeLowercase{\textit{et al.}}: Exact Analysis of Synchronizability for Complex Networks using Regular Graphs }
%



\maketitle

\begin{abstract}
Network synchronization is an emerging phenomenon in complex networks. The spectrum of Laplacian matrix will be immensely helpful for getting the network dynamics information.  Especially, network synchronizability is characterized by the ratio of second smallest eigen value to largest eigen value of the Laplacian matrix. We study the synchronization of complex networks modeled by regular graphs. We obtained the analytical expressions for network synchronizability for r-nearest neighbor cycle and r-nearest neighbor torus. We have also derived the generalized expression for synchronizability for m-dimensional r-nearest neighbor torus. The obtained analytical results agree with the simulation results and shown the effect of network dimension, number of nodes and overhead on synchronizability and connectivity in complex networks. This work provides the basic analytical tools for managing and controlling the synchronization in the finite sized complex networks and also given the generalized expressions for eigen values of Laplacian matrix for multi dimensional $r$- nearest neighbor networks.  
\end{abstract}

\begin{IEEEkeywords}
Complex networks, Syncronizability, Laplacian Spectra, Regular graphs
\end{IEEEkeywords}

%
\IEEEpeerreviewmaketitle

\section{Introduction}

\IEEEPARstart{S}{tudying} and understanding the synchronization processes is the intense research area in physical, biological, chemical, technological and social systems \cite{structure}.  Network synchronizability has been widely studied in the literature and still it is an exciting research problem \cite{stro}. In \cite{guan}, \cite{chen} generalized synchronization has been studied for complex networks. 

Investigating the network parameters which effect the synchronization process gives the important insights to understanding the network dynamics \cite{factors}. Synchronizability is only determined by the average degree and does not depend the system size and degree distribution \cite{motter}, whereas, maximum synchronizability of an network is completely determined by its associated feedback system \cite{jlu}. In \cite{lu}, it is proved that dynamic behaviors at each node and coupling configuration influence the network synchronization, also the synchronized region is directly related to inner linking matrix and synchronization is easy for larger values of algebraic connectivity  and ratio of eigen values \cite{chen}. The topology features of complex network greatly influence the dynamics processes of the network \cite{structure}. The work in \cite{bara} proved that small world networks synchronize as efficiently as random graphs and hyper cubes. 

In our work we prove that, network synchronization is greatly effected by the number of nodes and overhead (nearest neighbors) and network dimension. The Laplacian spectrum of a network plays a vital role to study the network dynamics towards synchronization \cite{spectra}. We derived the generalized expressions for eigen values of Laplacian, which can also be used to study the network dynamics. By using the algebraic graph theory, network synchronizability of large class of networks can be characterized by the ratio of second smallest eigen value to largest eigen value of the laplacian matrix \cite{bara}, \cite{control}, \cite{cheng}.  

We analyzed the synchronization of complex networks modeled by regular graphs with varying number of nearest neighbors.  In r-nearest neighbor cycle and torus, an edge will be existed between every pair of neighbors that are r hops away. Variable $r$ in the $r$-nearest neighbor networks captures the levels of overhead in the network operations \cite{analysis}. In this paper, we derive a general formula to efficiently and exactly compute the Network synchronizability and connectivity for $r$- nearest neighbor cycle, $r$- nearest neighbor torus and $m$- dimesional torus networks. The advantage of this kind of analysis is the producing the precise results without depending on thousands of simulation trails \cite{analysis}. We model the WSN as regular graphs and derived the analytical formulas for network connectivity and synchronizability.\\
In summary, this paper is organized as follows. \\
1) In Section II, we have given brief review about the network synchronizability. \\
2) In Section III, we have discussed the distance regular graphs and derived the expressions for eigen values of Laplacain matrix for the $r$- nearest neighbor networks. \\
3) In Section IV, exact formulas for network synchronizability and connectivity have been derived. \\
4) In Section V, simulation results have been presented and compared with the obtained analytical results. 

\section{Synchronization in Networks}
Let $G= (V,E)$, be an undirected graph with node set $V = \left\{ {1,2,......n} \right\}$ and an edge set $E \subseteq V \times V$. Further, let $A$ be $n\times n$ symmetric adjacency matrix of the graph $G$, each entry of adjacency matrix is represented by $a_{ij}$, which is $1$ if node $i$ is connected to node $j$, else it is $0$.

The degree matrix $D$ is defined as the diagonal matrix whose entry is $d_{ii}  = \deg (v_i )$. If we consider a network consisting of $N$ nodes, where $x_i$ represents $m$-dimensional vector for the $i$-th node. The dynamics of every node can be described by 

\begin{equation}
\dot x = f(x) - \sigma L \otimes Hx
\label{1}
\end{equation}

where $f$ represents the oscillator equations, $\sigma$ is the overall coupling strength, $H$ is the coupling matrix and $L$ is the Laplacian matrix describes the topology of the network.
The Laplacian matrix of the graph $G$ is the $n \times n$ symmetric matrix $L=D-A$, whose entries are

\begin{equation}
l_{ij}  = l_{ji}  = \left\{ \begin{array}{l}
 \deg (v_i )\,\,\,if\,\,j = i \\ 
  - a_{ij} \,\,\,\,\,\,\,\,\,if\,\,j \ne \,i \\ 
 \end{array} \right.
 \label{2}
\end{equation}

$Remark$ 1: Synchronizability of a large class of networks is determined by the eigen value ratio \cite{bara}, i.e.
\begin{equation}
R  = \frac{{\lambda _1 (L)}}{{\lambda _{n-1} (L)}}
\label{3}
\end{equation}
Let $\lambda _1 (L)$ and $\lambda _{n-1} (L)$ are the second smallest, second largest eigen values of the laplacian matrix. Second smallest eigen value of a Laplacian matrix $(\lambda _1 (L))$ defined as the algebraic connectivity of graphs \cite{fiedler}.

\section{Regular Graphs}
\subsection{$r$-nearest neighbor cycle}
The $r$-nearest neighbor cycle $C_n^r$ can be represented by a circulant matrix \cite{Circul}. A circulant matrix is defined as
 
\begin{equation}
\left[ \begin{array}{l}
 a_1 \,\,a_2 \,\,........a_{n - 1} \,\,a_n  \\ 
 a_n \,\,a_1 \,\,....\,....a_{n - 2} \,a_{n - 1}  \\ 
 .\,\,\,\,\,\,\,.\,\,\,\,\,\,\,\,\,\,\,\,\,\,\,\,\,\,\,\,\,\,\,\,.\,\,\,\,\,\,\,\,\,. \\ 
 .\,\,\,\,\,\,\,.\,\,\,\,\,\,\,\,\,\,\,\,\,\,\,\,\,\,\,\,\,\,\,\,.\,\,\,\,\,\,\,\,\,. \\ 
 a_3 \,\,a_4 \,\,\,...........a_1 \,\,\,a_2  \\ 
 a_2 \,\,a_3 \,\,.............a_n \,\,a_1 \,\, \\ 
 \end{array} \right]
 \label{4}
\end{equation}

and $j$-th eigen value  of a circulant matrix can be expressed as

\begin{equation}
\lambda _j  = a_1  + a_2 \omega ^{j}  + .............. + a_n \omega ^{(n - 1)j} 
\label{5}
\end{equation}

where $\omega$ be the $n$-th root of 1. Then $\omega$ is the complex number: 
\begin{equation}
\omega  = \cos \left( {\frac{{2\pi }}{n}} \right) + i\sin \left( {\frac{{2\pi }}{n}} \right)= e^\frac{i2\pi}{n}
\label{6}
\end{equation}

The 1-nearest cycle and 2-nearest cycle are shown in Fig.1 and Fig.2 respectively. Let the adjacency matrix $A$ and the degree matrix $D$ of 1-nearest cycle, then they can be written as

\begin{equation}
A=\left[ \begin{array}{l}
 0\,\,1\,\,0\,\,..............0\,\,1 \\ 
 1\,\,0\,\,1\,\,..............0\,0 \\ 
 .\,\,\,\,.\,\,\,\,.\,\,\,\,\,\,\,\,\,\,\,\,\,\,\,\,\,\,\,\,\,\,\,\,\,.\,\,\,. \\ 
 .\,\,\,\,.\,\,\,\,.\,\,\,\,\,\,\,\,\,\,\,\,\,\,\,\,\,\,\,\,\,\,\,\,\,.\,\,\,. \\ 
 0\,\,0\,\,0\,\,..............0\,\,1 \\ 
 1\,\,0\,\,0\,\,..............1\,0\, \\ 
 \end{array} \right]
 \label{7}
\end{equation}

\begin{equation}
D=\left[ \begin{array}{l}
 2\,\,0\,0\,\,............0\,\,0 \\ 
 0\,\,2\,0\,\,....\,........0\,\,0 \\ 
 .\,\,\,\,.\,\,\,\,\,\,\,\,\,\,\,\,\,\,\,\,\,\,\,\,\,\,\,\,\,\,\,\,\,.\,\,\,. \\ 
 .\,\,\,\,.\,\,\,\,\,\,\,\,\,\,\,\,\,\,\,\,\,\,\,\,\,\,\,\,\,\,\,\,\,.\,\,\,. \\ 
 0\,\,0\,\,\,...............2\,0\,\,\, \\ 
 0\,\,0\,\,\,...............0\,\,2 \\ 
 \end{array} \right]
 \label{8}
\end{equation}

$Theorem$ 1:The generalized expression for eigenvalues of Laplacian matrix for $1$-nearest neighbor cycle $C_n^1$ can be expressed as,
\begin{equation}
\lambda _j (L(C_n^1))  = 2 - 2\cos \left( {\frac{{2\pi j}}{n}} \right)
\label{9}
\end{equation}
where $j=0,1,...(n-1)$.\\
$Proof$: From (\ref{7}) and (\ref{8}), the Laplacian matrix $L$ for $C_n^1 $ can be written as,
\begin{equation}
L=\left[ \begin{array}{l}
 \,\,2\,\,\,\,-1\,\,\,\,\,\,\,\,\,\,\,0\,\,.............0\,\,-1 \\ 
 -1\,\,\,\,\,\,\,\,2\,\,\,\,-1\,\,.............0\,\,\,\,\,\,\,\,\,0 \\ 
 \,\,\,.\,\,\,\,\,\,\,\,\,\,\,\,\,.\,\,\,\,\,\,\,\,\,\,\,\,.\,\,\,\,\,\,\,\,\,\,\,\,\,\,\,\,\,\,\,\,\,\,\,\,\,.\,\,\,\,\,\,\,\,\,\,. \\ 
\,\,\, .\,\,\,\,\,\,\,\,\,\,\,\,\,.\,\,\,\,\,\,\,\,\,\,\,\,.\,\,\,\,\,\,\,\,\,\,\,\,\,\,\,\,\,\,\,\,\,\,\,\,\,.\,\,\,\,\,\,\,\,\,\,. \\ 
\,\,\, 0\,\,\,\,\,\,\,\,\,\,\,0\,\,\,\,\,\,\,\,\,\,\,0\,\,.............2\,-1 \\ 
 -1\,\,\,\,\,\,\,\,\,0\,\,\,\,\,\,\,\,\,\,\,0\,\,.........-1\,\,\,\,\,\,\,\,2\, \\ 
 \end{array} \right]
 \label{10}
\end{equation}

Hence, by using (\ref{5}) and (\ref{10}), $j^{th}$ eigen value of Laplacian matrix for $C_n^1$ can be written as,

\begin{equation}
\lambda _j (L(C_n^1)) = 2 - 2\cos \left( {\frac{{2\pi j}}{n}} \right)\nonumber
\end{equation}
where $j=0,1,...(n-1)$.\\

$Theorem$ 2: The generalized expression for eigenvalues of Laplacian matrix $L$ for $r$-nearest neighbor cycle $C_n^r $ can be expressed as,
\begin{equation}
\lambda _j (L(C_n^r))\, = 2r + 1 - \frac{{\sin \frac{{(2r + 1)\pi j}}{n}}}{{\sin \frac{{\pi j}}{n}}}
\label{11}
\end{equation}
where $j=0,1,...(n-1)$.\\
$Proof$:  From (\ref{9}), we can observe that, the first row is enough to obtain the eigen values of any circulant matrix.

The first row of adjacency matrix (A), degree matrix (D) and Laplacian matrix (L) can be written as follows,
\begin{equation}
A_{1n}  = \left[ {0\underbrace {\,1\,\,1\,\,1\,........0..........\,1\,\,1\,\,1}_{2r\,times}\,} \right]
\label{12}
\end{equation} 

\begin{equation}
D_{1n}  = \left[ {2r\,\,0\,\,0\,\,0\,.....0\,\,0\,\,0\,\,0} \right]
\label{13}
\end{equation}

\begin{equation}
\resizebox{.9 \hsize} {!} {$L_{1n}  = \left[ {2r\,\,\underbrace { - 1\,\, - 1\,\, - 1\,.....0...\,.... - 1\, - 1\, - 1}_{2r\,times}} \right]$}
\label{14}
\end{equation}

By using (\ref{14}) and (\ref{5}), we can write the   

\begin{equation} 
\lambda _j (L(C_n^r)) = 2r - 2\sum\limits_{i = 1}^r {\cos \left( {\frac{{2\pi ji}}{n}} \right)}
\label{15}
\end{equation}

Lemma 1: Trigonometric identity of Dirichlet kernel \cite{ident}
\begin{equation}
1 + 2\sum\limits_{j = 1}^r {\cos (jx)}  = \frac{{\sin \left( {r + \frac{1}{2}} \right)x}}{{\sin \left( {\frac{x}{2}} \right)}}
\label{16}
\end{equation}

Hence, from the Lemma 1, (\ref{15}) can be rewritten as,

\begin{equation}
\lambda _j (L(C_n^r))\, = 2r + 1 - \frac{{\sin \frac{{(2r + 1)\pi j}}{n}}}{{\sin \frac{{\pi j}}{n}}} \nonumber
\end{equation}

\begin{figure}[!t]
\centering
\includegraphics[width=2.2in]{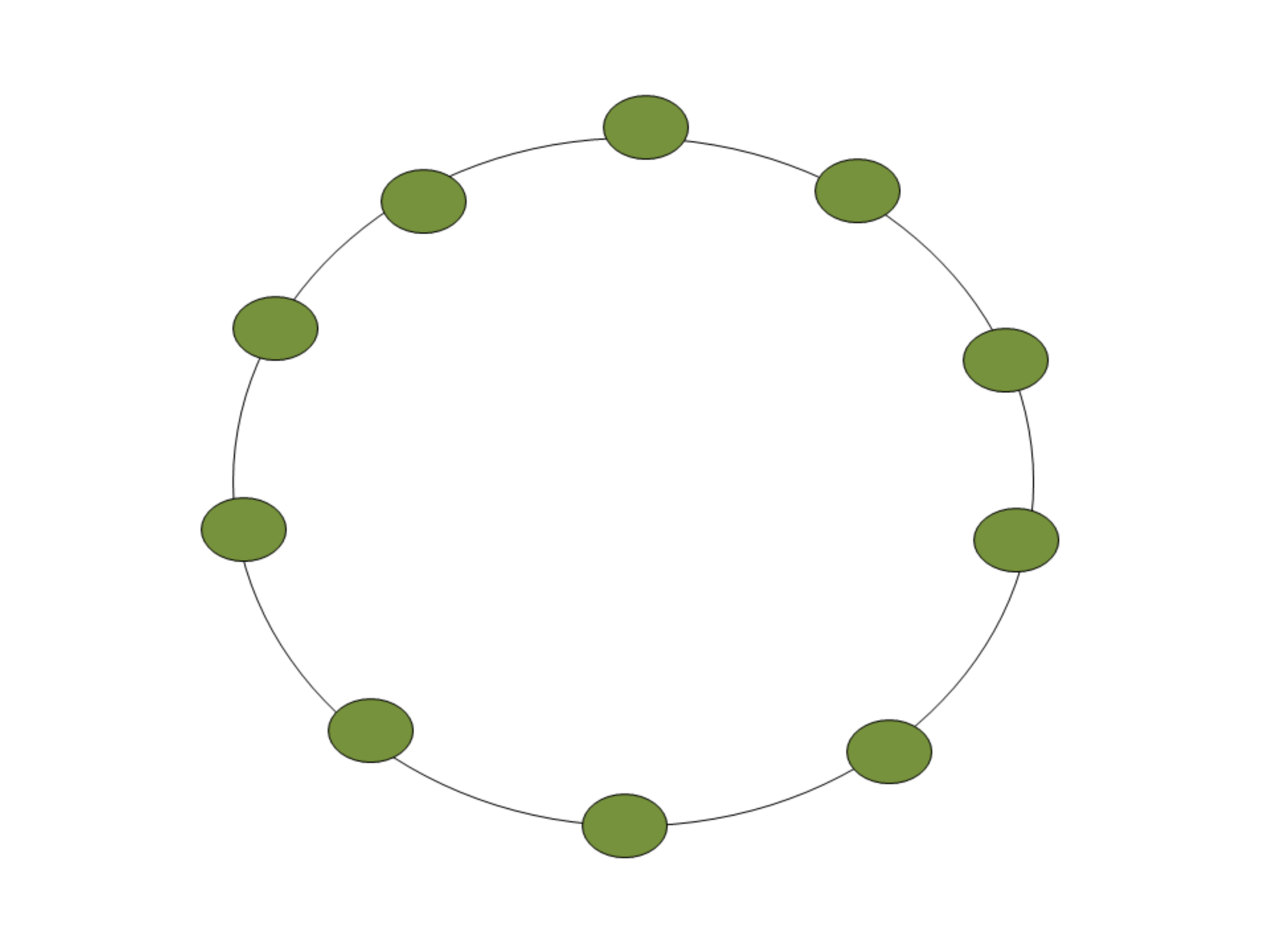}
\caption{1-nearest neighbor cycle }
\label{fig:1}
\end{figure}

\begin{figure}[!t]
\centering
\includegraphics[width=2in]{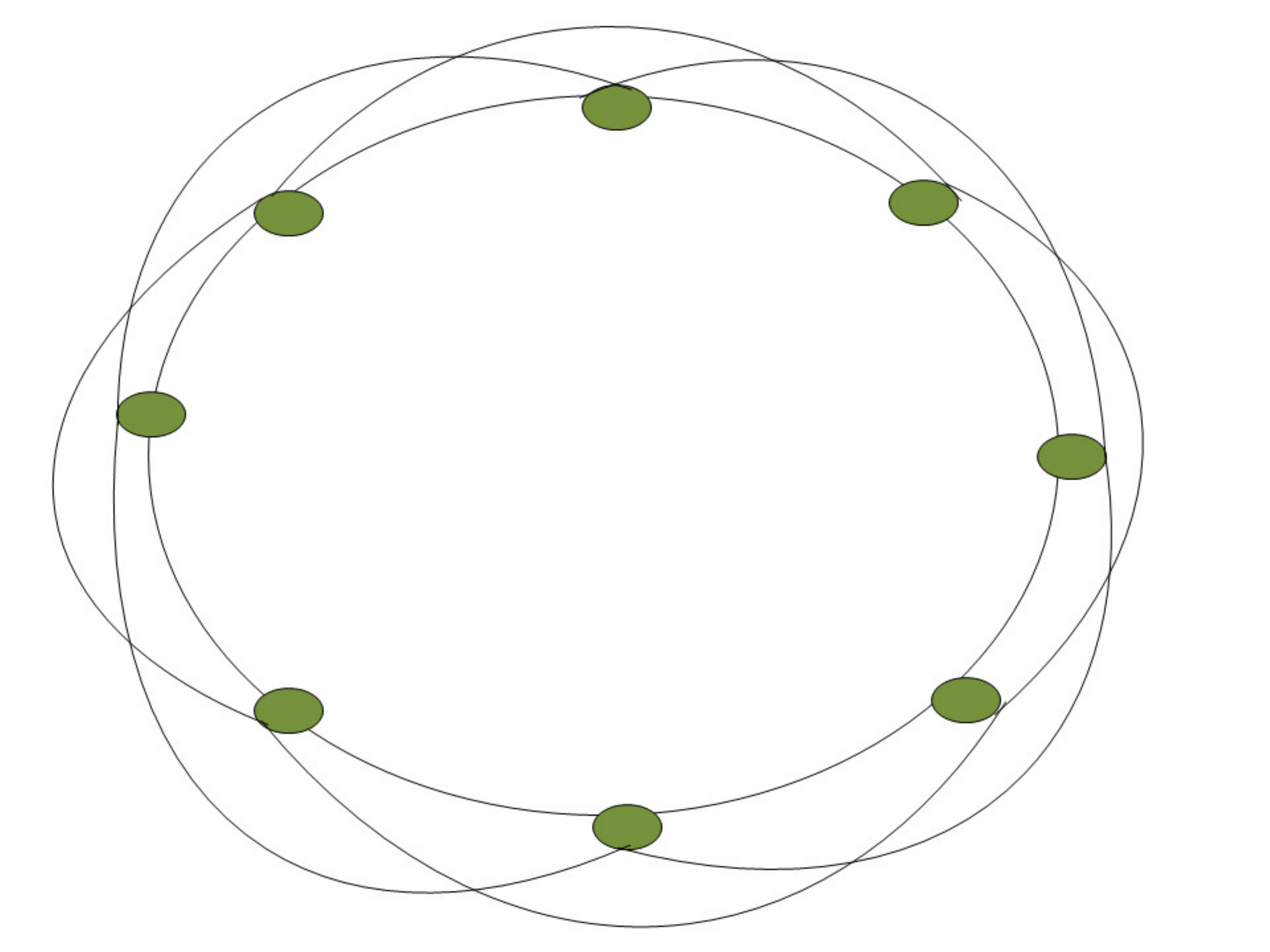}
\caption{2-nearest neighbor cycle}
\label{fig:2}
\end{figure}

\subsection{$r$-nearest neighbor torus}
A torus can be seen in Fig. 3 and it can be represented by the $n \times n$ block circulant matrix $A$ as
\begin{equation}
A = \left[ \begin{array}{l}
 A_0 \,\,\,\,\,\,\,\,\,\,A_1 \,\,........A_{n_1 - 2} \,\,A_{n_1 -1}  \\ 
 A_{n_1 -1} \,\,A_0 \,\,....\,....A_{n_1 - 3} \,A_{n_1 - 2}  \\ 
 .\,\,\,\,\,\,\,\,\,\,\,\,\,\,\,\,\,\,\,\,\,.\,\,\,\,\,\,\,\,\,\,\,\,\,\,\,\,\,\,\,\,\,\,\,\,.\,\,\,\,\,\,\,\,\,\,\,\,\,\,\,. \\ 
 .\,\,\,\,\,\,\,\,\,\,\,\,\,\,\,\,\,\,\,\,\,.\,\,\,\,\,\,\,\,\,\,\,\,\,\,\,\,\,\,\,\,\,\,\,\,.\,\,\,\,\,\,\,\,\,\,\,\,\,\,\,. \\ 
 A_1 \,\,\,\,\,\,\,\,\,\,\,\,\,A_2 \,\,..........A_{n_1 -1} \,\,A_0 \,\, \\ 
 \end{array} \right]
 \label{19}
\end{equation}

where the number of nodes $n=n_1^2$, then each block $A_i$, for $i=0,1...(n_1-1)$ represents $n_1 \times n_1$ circulant matrices.\\

$Lemma$ 2: Let $G$ be the cartesian product of two graphs $G^{'}$ and $G^{''}$ with vertex sets $V^{'}$ and $V^{''}$ and edge sets $E^{'}$ and $E^{''}$.
Let the eigen values of $G^{'}$ are  $\lambda _1 \left( {G^{'} } \right)..........\lambda _p \left( {G^{'}} \right)$ and $G^{''}$ are $\lambda _1 \left( {G^{''} } \right)..........\lambda _q \left( {G^{''} } \right)$, where $ p = \left| {V^{'} } \right|$ and $q = \left| {V^{''} } \right|$. Let the vertex set of $G$ is $r = \left| {V } \right|$, which can be expressed as $V = \left| {V^{'} } \right| \times \left| {V^{''} } \right|$ \cite{spec}. Then, the eigen values of $G$ can be expressed as 
\begin{equation}
\lambda _k \left( G \right) = \lambda _i \left( {G^{'} } \right) + \lambda _j \left( {G^{''} } \right)
\label{20}
\end{equation}
,where $i \in \{ 1,2,....p\}$, $j \in \{ 1,2,....q\}$ and $k \in \{ 1,2,....r\}$. \\

$Remark$ 2: (\ref{20}) also holds for eigen values of the Laplacians $L^{'}$ and $L^{''}$ of graphs of $G^{'}$ and $G^{''}$ respectively\cite{opt}.\\

$Theorem$ 3: The generalized expression for eigenvalues of Laplacian matrix $L$ for $1$-nearest neighbor torus $T_n^1 $ can be expressed as

\begin{equation}
\lambda _{j_1 ,j_2 } \left( {L(T_{k_1 ,k_2 }^1 )} \right) = 4 - 2\cos \left( {\frac{{2\pi j_1 }}{{k_1 }}} \right) - 2\cos \left( {\frac{{2\pi j_2 }}{{k_2 }}} \right)
\label{21}
\end{equation}
where $j_1  = 0,1,2,...(k_1  - 1), j_2  = 0,1,2,...(k_2  - 1)$.\\

$Proof$ : $1$-nearest neighbor two dimensional torus can be represented by the Cartesian product of two $1$-nearest neighbor cycles.
So from the $Lemma$ 2, we can write the $\lambda _{j_1 ,j_2 } \left( {L(T_{k_1 ,k_2 }^1 )} \right)$ as,

\begin{equation}
\lambda _{j_1 ,j_2 } \left( {L(T_{k_1 ,k_2 }^1 )} \right) = \lambda _{j_1 } \left( {L(C_{k_1 }^1 )} \right) + \lambda _{j_2 } \left( {L(C_{k_2 }^1 )} \right)
\label{22}
\end{equation}

To prove the theorem, write the expressions for $\lambda _{j_1 } \left( {L(C_{k_1 }^1 )} \right)$ and $\lambda _{j_2 } \left( {L(C_{k_2 }^1 )} \right)$ using (\ref{9}) and substitute in (\ref{22}). \\

$Theorem$ 4: The generalized expression for eigenvalues of Laplacian matrix $L$ for $r$-nearest neighbor torus $T_n^r$ can be expressed as
\begin{equation}
\lambda _{j_1 ,j_2 } \left( {L(T_{k_1 ,k_2 }^1 )} \right) = 4r + 2 - \frac{{\sin \frac{{(2r + 1)\pi j_1 }}{{k_1 }}}}{{\sin \frac{{\pi j_1 }}{{k_1 }}}} - \frac{{\sin \frac{{(2r + 1)\pi j_2 }}{{k_2 }}}}{{\sin \frac{{\pi j_2 }}{{k_2 }}}}
\label{23}
\end{equation}

where $j_1  = 0,1,2,...(k_1  - 1), j_2  = 0,1,2,...(k_2  - 1)$.\\

$Proof$ :  
$T_n^r$ can be represented by Cartesian product of two $r$-nearest neighbor cycles.
So from the $Lemma$ 2, we can write the $\lambda _{j_1 ,j_2 } \left( {L(T_{k_1 ,k_2 }^r )} \right)$ as,

\begin{equation}
\lambda _{j_1 ,j_2 } \left( {L(T_{k_1 ,k_2 }^r )} \right) = \lambda _{j_1 } \left( {L(C_{k_1 }^r )} \right) + \lambda _{j_2 } \left( {L(C_{k_2 }^r )} \right)
\label{24}
\end{equation}

From (\ref{11}), we can write the expressions for $\lambda _{j_1 } \left( {L(C_{k_1 }^r )} \right)$ and  $\lambda _{j_2 } \left( {L(C_{k_2 }^r )} \right)$, substituting them in (\ref{24}) proves the theorem.
\\

$Theorem$ 5: The generalized expression for eigenvalues of Laplacian matrix $L$ for $m$-dimensional $r$-nearest neighbor torus can be expressed as

\begin{equation}
\resizebox{.9 \hsize} {!} {$\lambda _{j_1 ,j_2,... j_m} \left( {L(T_{k_1 ,k_2....k_m}^r )} \right) = (2r + 1)m - \sum\limits_{i = 1}^m {\left( {\frac{{\sin \frac{{(2r + 1)\pi j_i }}{{k_i }}}}{{\sin \frac{{\pi j_i }}{{k_i }}}}} \right)}$}
\label{25}
\end{equation}

$Proof$ 5: 
$r$-nearest neighbor $m$-dimensional torus can be represented by Cartesian product of $m$ number of $r$-nearest neighbor cycles.
So from the $Lemma$ 2, we can write the $\lambda _{j_1 ,j_2....j_m } \left( {L(T_{k_1 ,k_2....k_m }^r )} \right)$ as,

\begin{equation}
\resizebox{.9 \hsize} {!} {$\lambda _{j_1 ,j_2,... j_m} \left( {L(T_{k_1 ,k_2....k_m}^r )} \right) = \lambda _{j_1 } \left( {L(C_{k_1 }^r )} \right) +\lambda _{j_2 } \left( {L(C_{k_1 }^r )} \right) .............+\lambda _{j_m } \left( {L(C_{k_m }^r )} \right)$}
\label{26}
\end{equation}

From (\ref{11}), we can substitute the expressions for $\lambda _{j_1 } \left( {L(C_{k_1 }^r )} \right)$, $\lambda _{j_2 } \left( {L(C_{k_2 }^r )} \right)$ and $\lambda _{j_m } \left( {L(C_{k_m }^r )} \right)$ in (\ref{26}), which proves the theorem.\\
\\

\begin{figure}[!t]
\centering
\includegraphics[width=2.2 in]{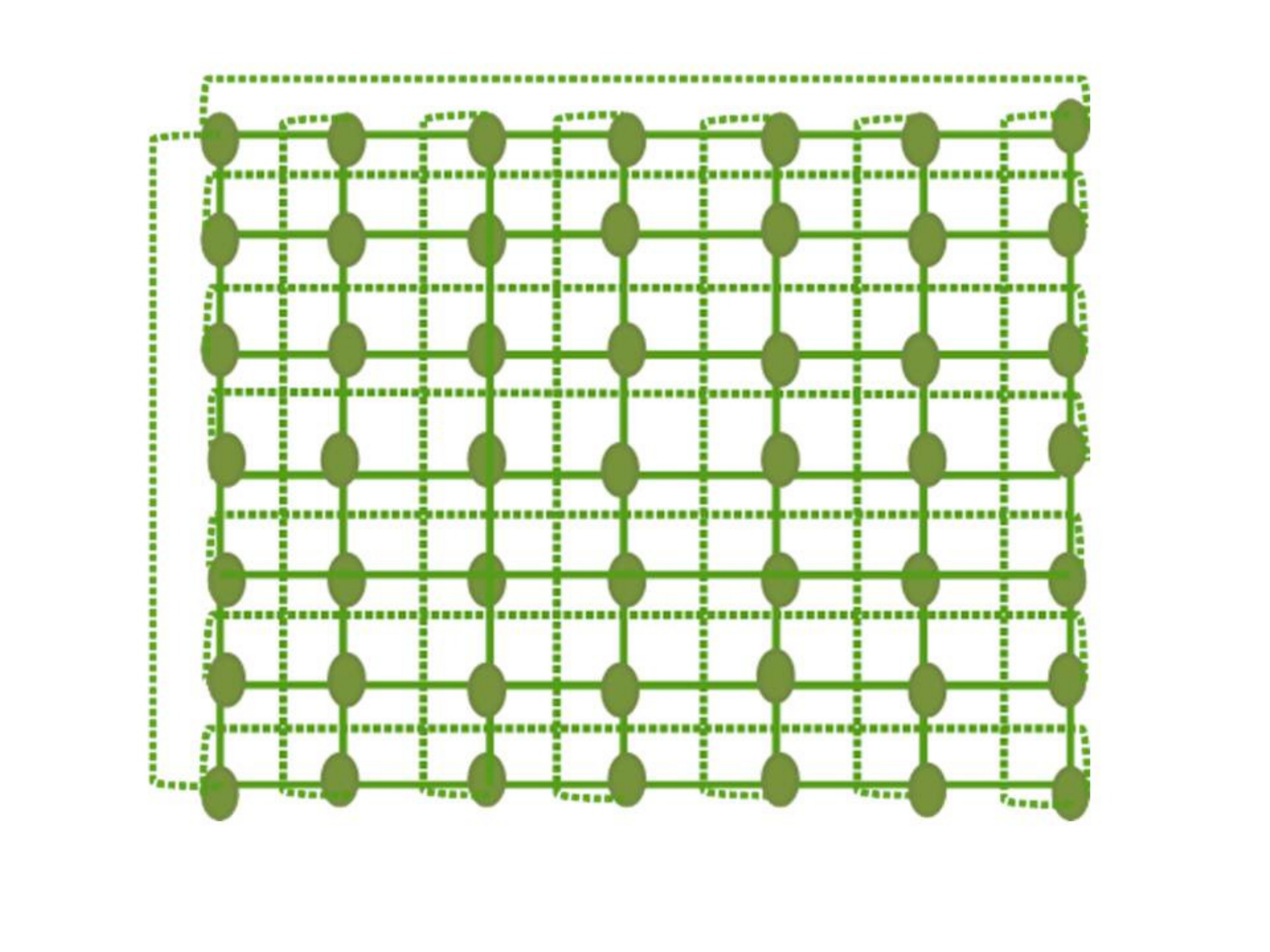}
\caption{Two dimensional torus }
\label{fig:3}
\end{figure}

\section{Synchronization analysis for regular graphs}

\subsection{$r$-nearest neighbor cycle}
$Theorem$ 6: Given $r$-nearest neighbor cycle $C_n^r $ and $n$ is even integer, then network synchronizability $R$ can be expressed as

\begin{equation}
R = \frac{{(2r + 1)sin\frac{\pi }{n} - \sin \frac{{(2r + 1)\pi }}{n}}}{{(2r + 1 - \cos \pi r)\sin \frac{\pi }{n}}}
\label{27}
\end{equation}

$Proof$ :
From (\ref{3}), network synchronizability can be expressed as the ratio of second smallest eigen value of Laplacian matrix to largest eigen value of Laplacian matrix . \\ Since $\lambda _{n/2} (L)$ is the largest eigen value of Laplacian matrix for $r$-nearest neighbor cycle $C_n^r $, synchronizability  can be rewritten as 

\begin{equation}
R=\frac{{\lambda _1 (L)}}{{\lambda _{n/2} (L)}}
\label{28}
\end{equation}

By substituting $j=1$ in (\ref{11}), we can write the $\lambda _1 (L)$ as,

\begin{equation}
\lambda _1 (L)\, = 2r + 1 - \frac{{\sin \frac{{(2r + 1)\pi }}{n}}}{{\sin \frac{\pi }{n}}}
\label{29}
\end{equation}

By substituting $j=n/2$ in (\ref{11}), we can write the $\lambda _{n/2} (L)$ as,
\begin{equation}
\lambda _{\frac{n}{2}} (L)\, = 2r + 1 - \cos \pi r
\label{30}
\end{equation}

Hence theorem is proved by substituting the (\ref{29}) and (\ref{30})in (\ref{28}). \\

$Theorem$ 7: Given $r$-nearest neighbor cycle $C_n^r $ and $n$ is odd integer, then network synchronizability $R$ can be expressed as 
\begin{equation}
R = \frac{{\left( {(2r + 1)sin\frac{\pi }{n} - \sin \frac{{(2r + 1)\pi }}{n}} \right)\cos \frac{\pi }{{2n}}}}{{\left( {(2r + 1)cos\frac{\pi }{{2n}} - \cos \frac{{(2r + 1)\pi }}{n}} \right)\sin \frac{\pi }{n}}}
\label{31}
\end{equation}

$Proof$ : 

$\lambda _{(n-1)/2} (L)$ is the largest eigen value for $r$-nearest neighbor cycle $C_n^r $ and (\ref{3}) can be rewritten as 
 
 \begin{equation}
R=\frac{{\lambda _1 (L)}}{{\lambda _{(n-1)/2} (L)}}
\label{32}
\end{equation}

By substituting $j=(n-1)/2$ in (\ref{11}), we can write the $\lambda _{(n-1)/2} (L)$ as,

\begin{equation}
\lambda _{\frac{{n - 1}}{2}} (L) = 2r + 1 - \frac{{\sin \frac{{(2r + 1)\pi (n - 1)}}{{2n}}}}{{\sin \frac{{\pi (n - 1)}}{{2n}}}}
\label{33}
\end{equation}

Hence theorem is proved by substituting the (\ref{29}) and (\ref{33}) in (\ref{32}).

\subsection{$r$-nearest neighbor torus}
$Theorem$ 8: Given a $r$-nearest neighbor torus  $T_n^r$ for $k_1, k_2$ are even integers, then network synchronizability $R$ can be expressed as
\begin{equation}
R = \frac{{2r + 1 - \left( {\frac{{\sin \frac{{(2r + 1)\pi }}{{k_2 }}}}{{\sin \frac{\pi }{{k_2 }}}}} \right)}}{{4r + 2 - 2\cos \pi r}}
\label{34}
\end{equation}
$Proof$ : 
$\lambda _{0,1} (L)$ is the second smallest eigen value and $\lambda _{\frac{{k_1 }}{2},\frac{{k_2 }}{2}} (L)$ is the largest eigen value of Laplacian matix for $r$-nearest neighbor torus for $k_1, k_2$ are even integers and network synchronizability $R$  can be rewritten as,

\begin{equation}
R=\frac{\lambda _{0,1} (L)}{\lambda _{\frac{{k_1 }}{2},\frac{{k_2 }}{2}} (L)}
\label{35}
\end{equation}

By substituting $j_{1}=0$ and $j_{2}=1$ in (\ref{21}), we can get 
\begin{equation}
\lambda _{0,1} (L)\, = 2r + 1 - \frac{{\sin \frac{{(2r + 1)\pi }}{{k_2 }}}}{{\sin \frac{\pi }{{k_2 }}}}
\label{36}
\end{equation}
By substituting $j_{1}=k_1/2 $ and $j_{2}=k_2/2$ in (\ref{21}), we can get 
\begin{equation}
\lambda _{\frac{{k_1 }}{2},\frac{{k_2 }}{2}} (L)\, = 4r + 2 - 2\cos \pi r
\label{37}
\end{equation}

Hence theorem is proved by substituting the (\ref{36}) and (\ref{37}) in (\ref{35}). \\

$Theorem$ 9: Given a $r$-nearest neighbor torus $T_n^r$ for $k_1, k_2$ are odd integers, then network synchronizability $R$ can be expressed as
\begin{equation}
\resizebox{.9 \hsize} {!} {$R = \frac{{2r + 1 - \left( {\frac{{\sin \frac{{(2r + 1)\pi }}{{k_2 }}}}{{\sin \frac{\pi }{{k_2 }}}}} \right)}}{{4r + 2 - \left( {\frac{{\sin \frac{{(2r + 1)\pi (k_1  - 1)}}{{2k_1 }}}}{{\sin \frac{{\pi (k_1  - 1)}}{{2k_1 }}}}} \right) - \left( {\frac{{\sin \frac{{(2r + 1)\pi (k_2  - 1)}}{{2k_2 }}}}{{\sin \frac{{\pi (k_2  - 1)}}{{2k_2 }}}}} \right)}}$}
\label{38}
\end{equation}
 
$Proof$ : 

$ \lambda _{\frac{{(k_1  - 1)}}{2},\frac{{(k_2  - 1)}}{2}} (L)$ is the largest eigen value of the Laplacian matrix for $r$-nearest neighbor torus when $k_1, k_2$ are odd integers. Then expression for network synchronizability can be rewritten as

\begin{equation}
R = \frac{{\lambda _{0,1} }}{{\lambda _{\frac{{(k_1  - 1)}}{2},\frac{{(k_2  - 1)}}{2}} }}
\label{39}
\end{equation}

By substituting $j_{1}=(k_1-1)/2 $ and $j_{2}=(k_2-1)/2$ in (\ref{21}), we can get 

\begin{equation}
\resizebox{.9 \hsize} {!} {$\lambda _{\frac{{(k_1  - 1)}}{2},\frac{{(k_2  - 1)}}{2}} (L) = 4r + 2 - \frac{{\sin \frac{{(2r + 1)\pi (k_1  - 1)}}{{2k_1 }}}}{{\sin \frac{{\pi (k_1  - 1)}}{{2k_1 }}}} - \frac{{\sin \frac{{(2r + 1)\pi (k_2  - 1)}}{{2k_2 }}}}{{\sin \frac{{\pi (k_2  - 1)}}{{2k_2 }}}}$}
\label{40}
\end{equation}

Hence theorem is proved by substituting the (\ref{36}) and (\ref{40}) in (\ref{39}).

\subsection{$m$-dimensional torus}
$Theorem$ 10: Given a $m$-dimensional $r$-nearest neighbor torus for $k_1, k_2, k_3...k_m$ are even integers, then network synchronizability $R$ can be expressed as
\begin{equation}
R = \frac{{(2r + 1) - \left( {\frac{{\sin \frac{{(2r + 1)\pi j_1 }}{{k_1 }}}}{{\sin \frac{{\pi j_1 }}{{k_1 }}}}} \right)}}{{m\left( {2r + 1 - \cos (\pi r)} \right)}}
\label{41}
\end{equation}
$Proof$ : 

$\lambda _{1,0,.....0}(L)$ is the second smallest eigen value and $\lambda _{\frac{{k_1 }}{2},\frac{{k_2 }}{2},.....\frac{{k_m }}{2}} (L)$ is the largest eigen value of Laplacian matix for $m$-dimensional $r$-nearest neighbor torus for $k_1, k_2, k_3...k_m$ are even integers and network synchronizability $R$  can be rewritten as,

\begin{equation}
R=\frac{\lambda _{1,0,.....0}(L)}{\lambda _{\frac{{k_1 }}{2},\frac{{k_2 }}{2},.....\frac{{k_m }}{2}} (L)}
\label{42}
\end{equation}

By substituting the $j_{1}=1$ and $j_{2}=j_{3}......j_{n}=0$ in (\ref{23}), we can get 
\begin{equation}
\lambda _{1,0,.....0} (L) = (2r + 1) - \left( {\frac{{\sin \frac{{(2r + 1)\pi j_1 }}{{k_1 }}}}{{\sin \frac{{\pi j_1 }}{{k_1 }}}}} \right)
\label{43}
\end{equation}

By substituting $j_{1}=k_1/2 $, $j_{2}=k_2/2$.....$j_{m}=k_m/2$ in (\ref{23}), we can get 
\begin{equation}
\lambda _{\frac{{k_1 }}{2},\frac{{k_2 }}{2},.....\frac{{k_m }}{2}} (L) = m\left( {2r + 1 - \cos (\pi r)} \right)
\label{44}
\end{equation}

Hence theorem is proved by substituting the (\ref{43}) and (\ref{44}) in (\ref{42}). \\

$Theorem$ 11: Given a $m$-dimensional $r$-nearest neighbor torus and for $k_1, k_2, k_3...k_m$ are odd integers, then network synchronizability can be expressed as
\begin{equation}
R = \frac{{(2r + 1) - \left( {\frac{{\sin \frac{{(2r + 1)\pi j_1 }}{{k_1 }}}}{{\sin \frac{{\pi j_1 }}{{k_1 }}}}} \right)}}{{(2r + 1)m - \sum\limits_{i = 1}^m {\left( {\frac{{\sin \frac{{(2r + 1)\pi (k_i  - 1)}}{{k_i }}}}{{\sin \frac{{\pi (k_i  - 1)}}{{k_i }}}}} \right)} }}
\label{45}
\end{equation}
$Proof$ : 
$\lambda _{\frac{{(k_1  - 1)}}{2},\frac{{(k_2  - 1)}}{2},.....\frac{{(k_m  - 1)}}{2}} (L)$ is the largest eigen value for a $m$-dimensional $r$-nearest neighbor torus when $k_1, k_2, k_3...k_m$ are odd integers, then network synchronizability can be rewritten as

\begin{equation}
R=\frac{\lambda _{1,0,.....0} (L)}{\lambda _{\frac{{(k_1  - 1)}}{2},\frac{{(k_2  - 1)}}{2},.....\frac{{(k_m  - 1)}}{2}} (L)}
\label{46}
\end{equation}

By substituting $j_{1}=(k_1  - 1)/2 $, $j_{2}=(k_2  - 1)/2$.....$j_{m}=(k_m  - 1)/2$ in (\ref{23}), we can get

\begin{equation}
\resizebox{.9 \hsize} {!} {$\lambda _{\frac{{(k_1  - 1)}}{2},\frac{{(k_2  - 1)}}{2},.....\frac{{(k_m  - 1)}}{2}} (L) = (2r + 1)m - \sum\limits_{i = 1}^m {\left( {\frac{{\sin \frac{{(2r + 1)\pi (k_i  - 1)}}{{k_i }}}}{{\sin \frac{{\pi (k_i  - 1)}}{{k_i }}}}} \right)}$} 
\label{47}
\end{equation}

Hence theorem is proved by substituting the (\ref{42}) and (\ref{47}) in (\ref{46}).

\section{Simulation results}
\begin{figure}[!t]
\centering
\includegraphics[width=2.5 in]{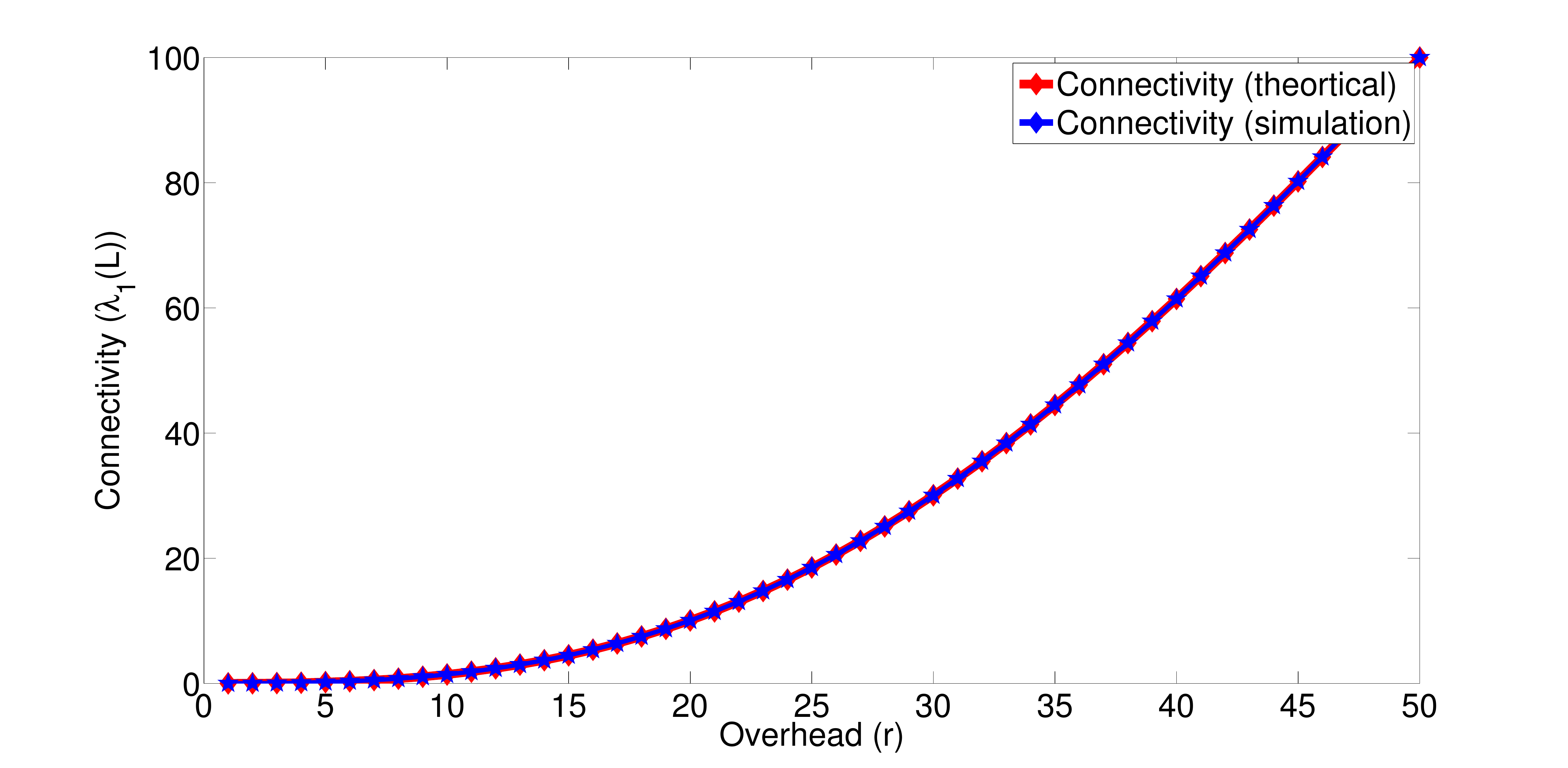}
\caption{Connectivity versus Overhead for $r$-nearest neighbor cycle}
\label{fig:5}
\end{figure}

\begin{figure}[!t]
\centering
\includegraphics[width=2.5 in]{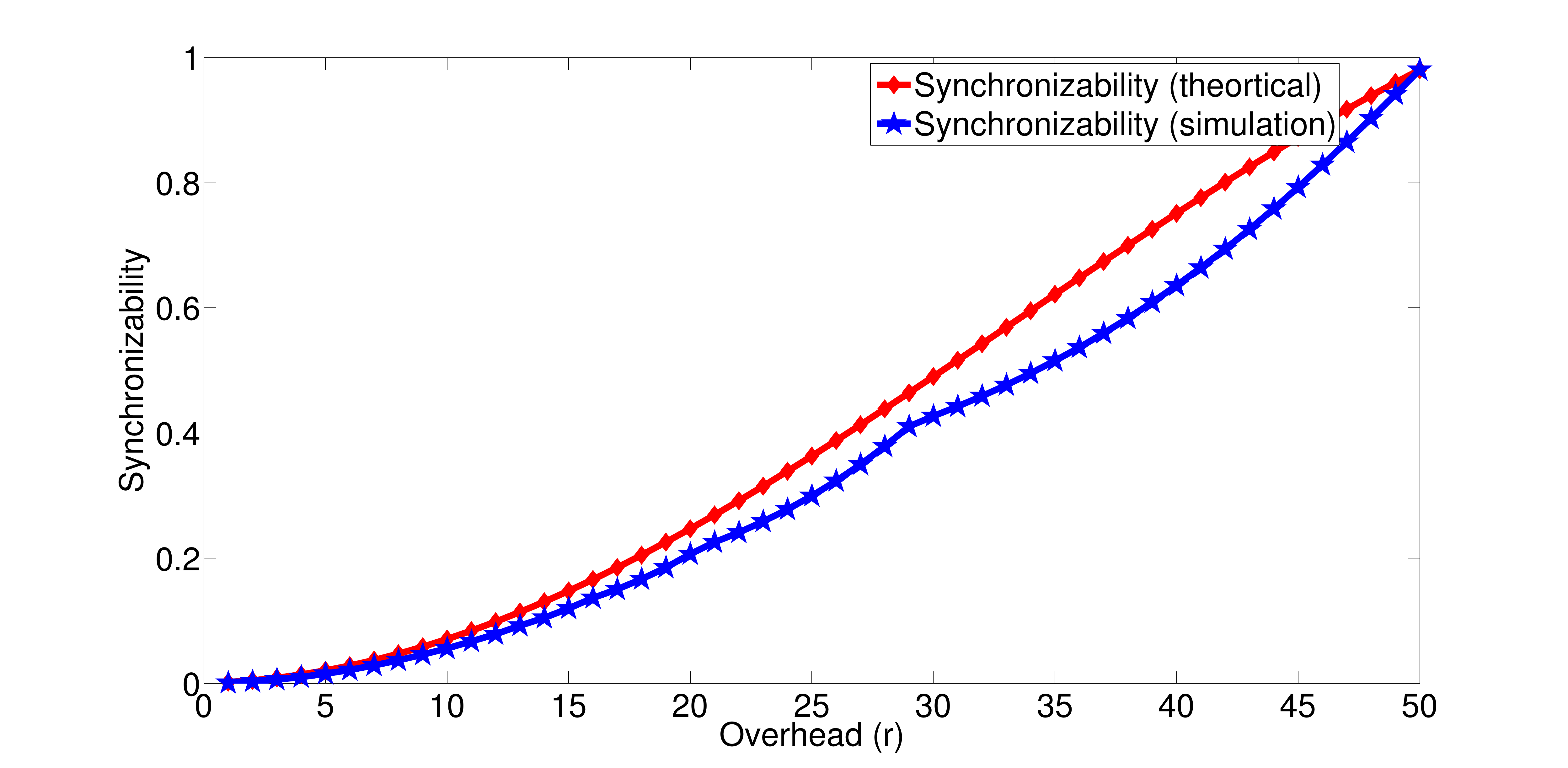}
\caption{Synchronizability versus Overhead for $r$-nearest neighbor cycle}
\label{fig:6}
\end{figure}

\begin{figure}[!t]
\centering
\includegraphics[width=2.5 in]{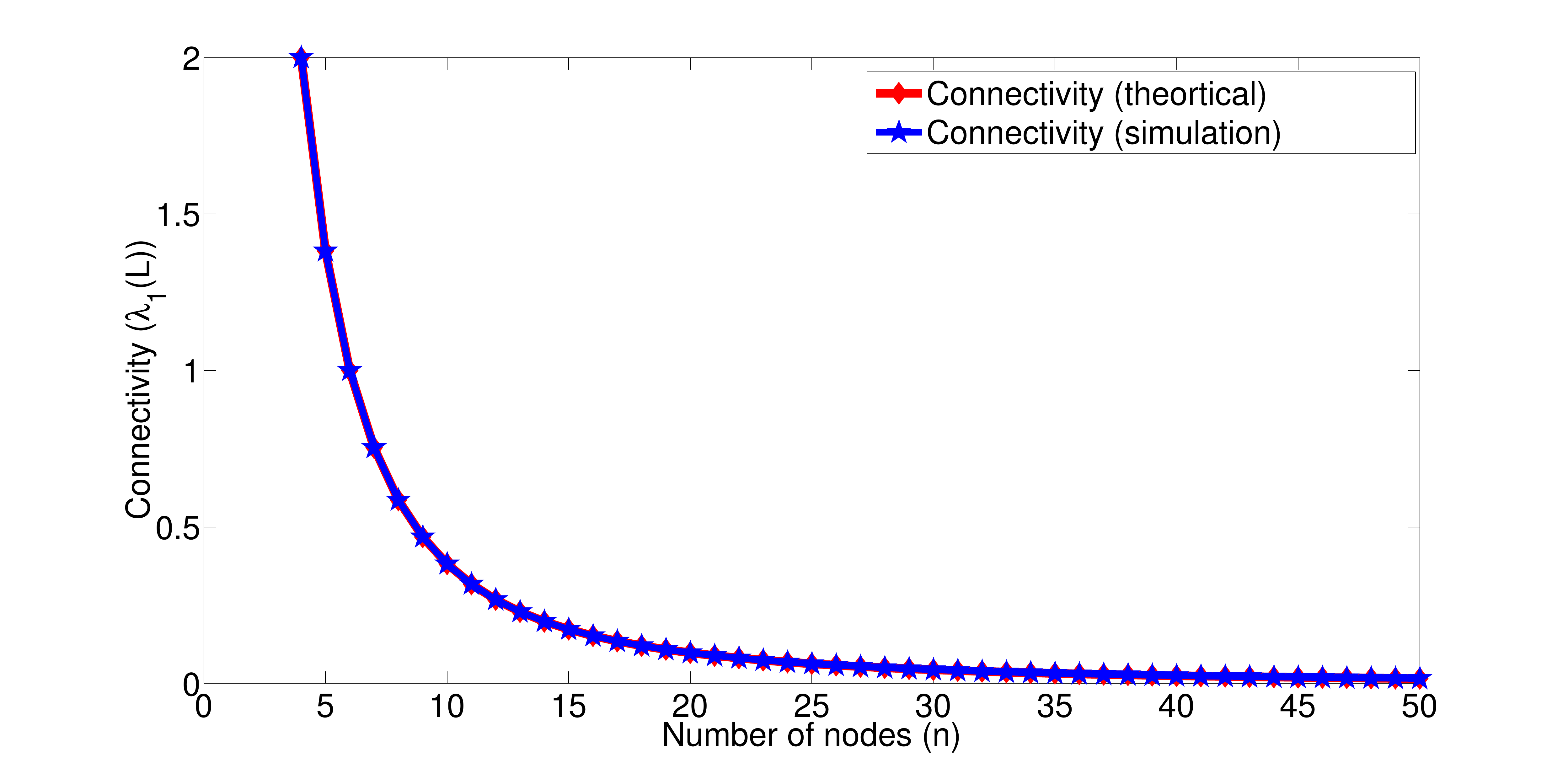}
\caption{Connectivity versus number of nodes for $r$-nearest neighbor cycle}
\label{fig:7}
\end{figure}

\begin{figure}[!t]
\centering
\includegraphics[width=2.5 in]{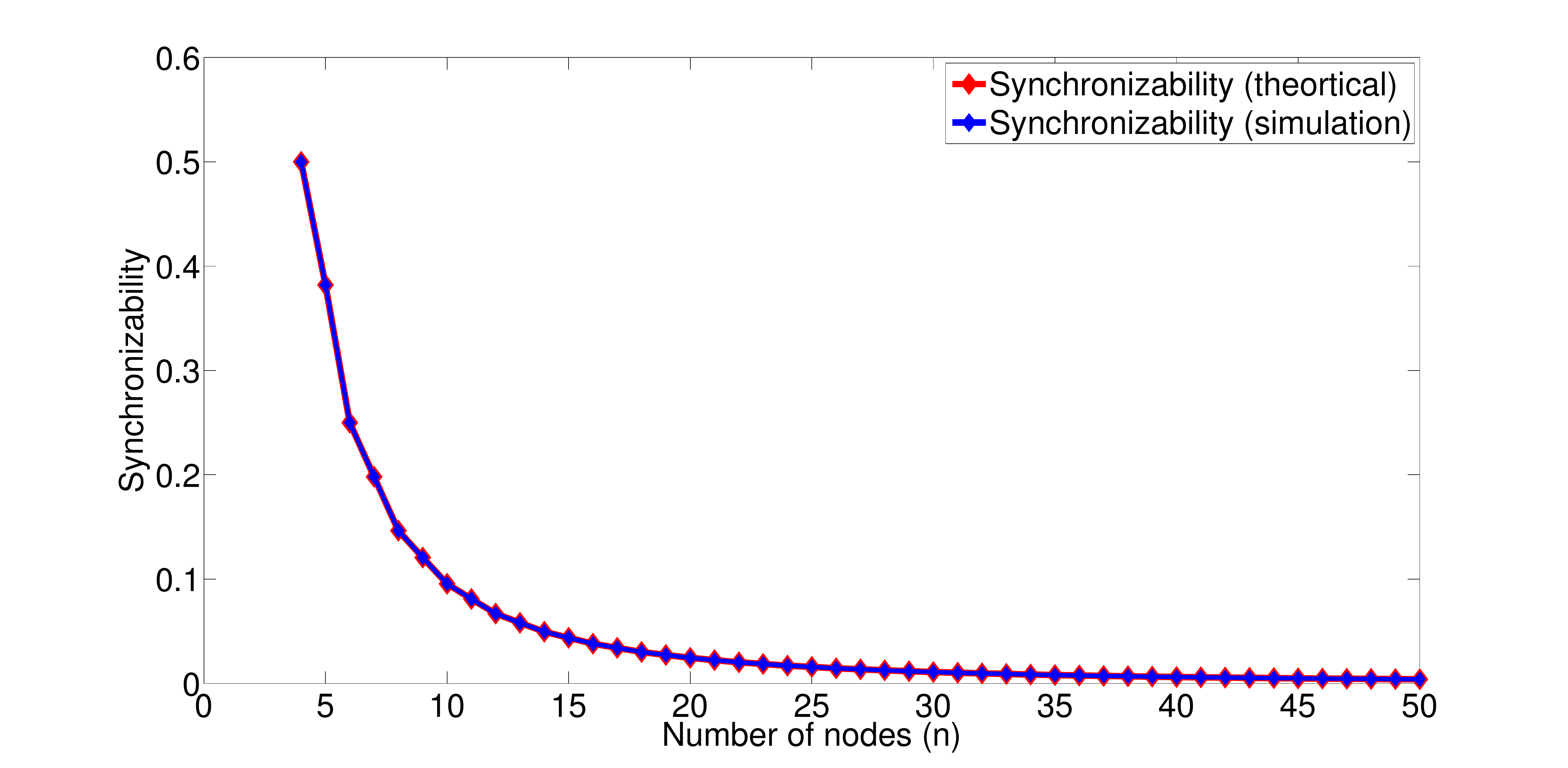}
\caption{Synchronizability versus number of nodes for $r$-nearest neighbor cycle}
\label{fig:8}
\end{figure}

\begin{figure}[!t]
\centering
\includegraphics[width=2.5 in]{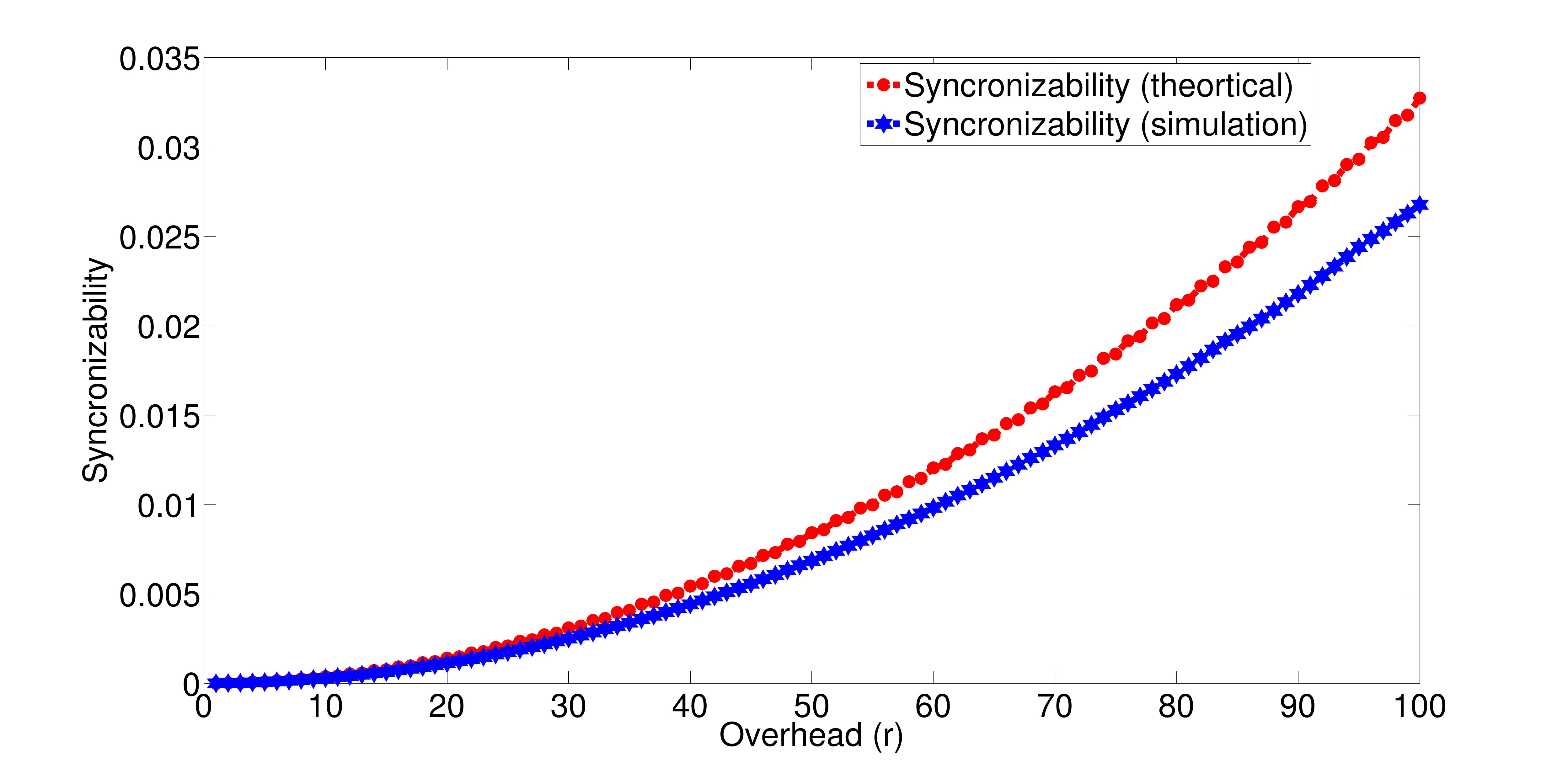}
\caption{Synchronizability versus Overhead for $r$-nearest neighbor torus}
\label{fig:9}
\end{figure}

\begin{figure}[!t]
\centering
\includegraphics[width=2.5 in]{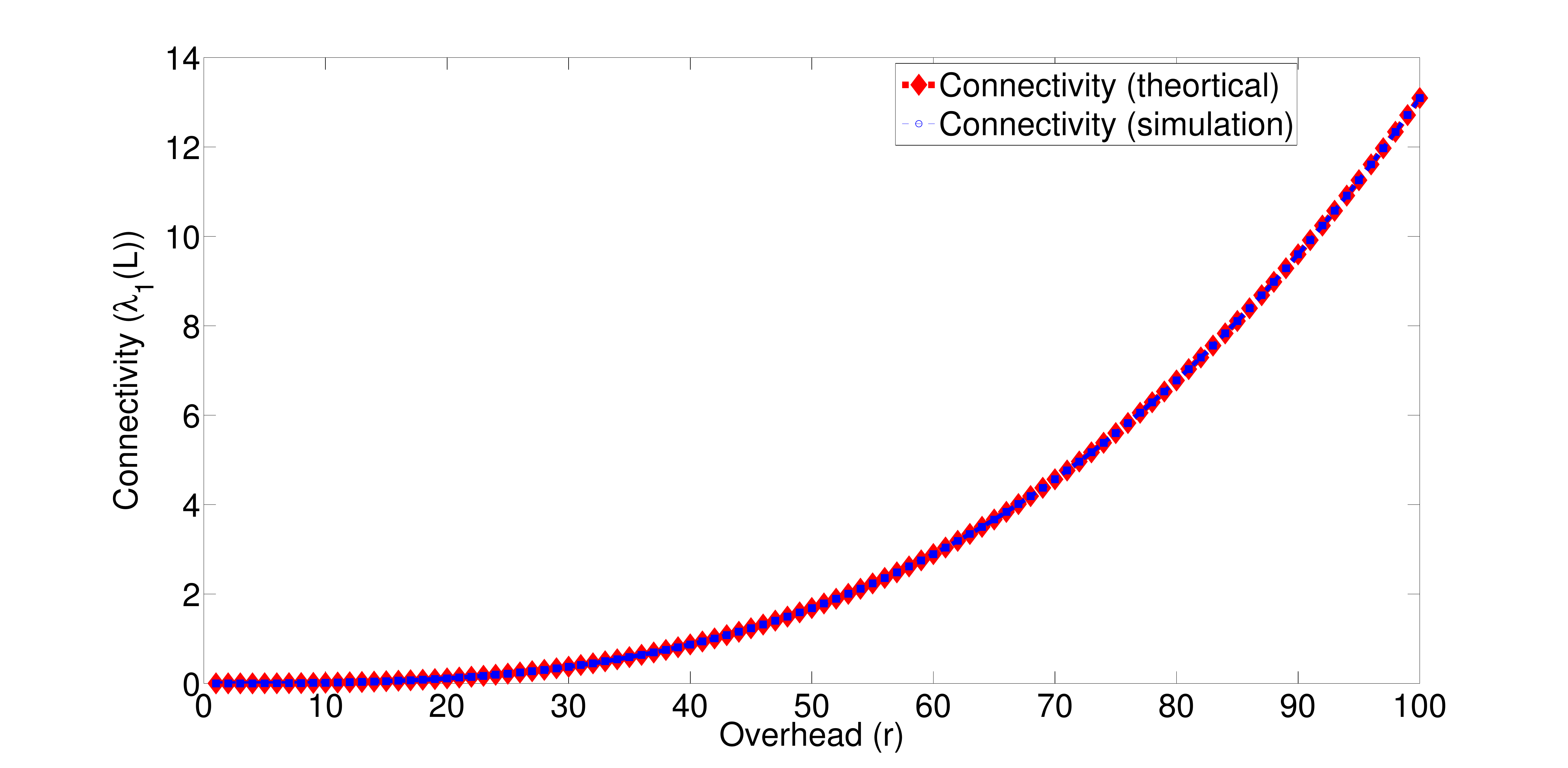}
\caption{Connectivity versus Overhead for $r$-nearest neighbor torus}
\label{fig:10}
\end{figure}

\begin{figure}[!t]
\centering
\includegraphics[width=2.5 in]{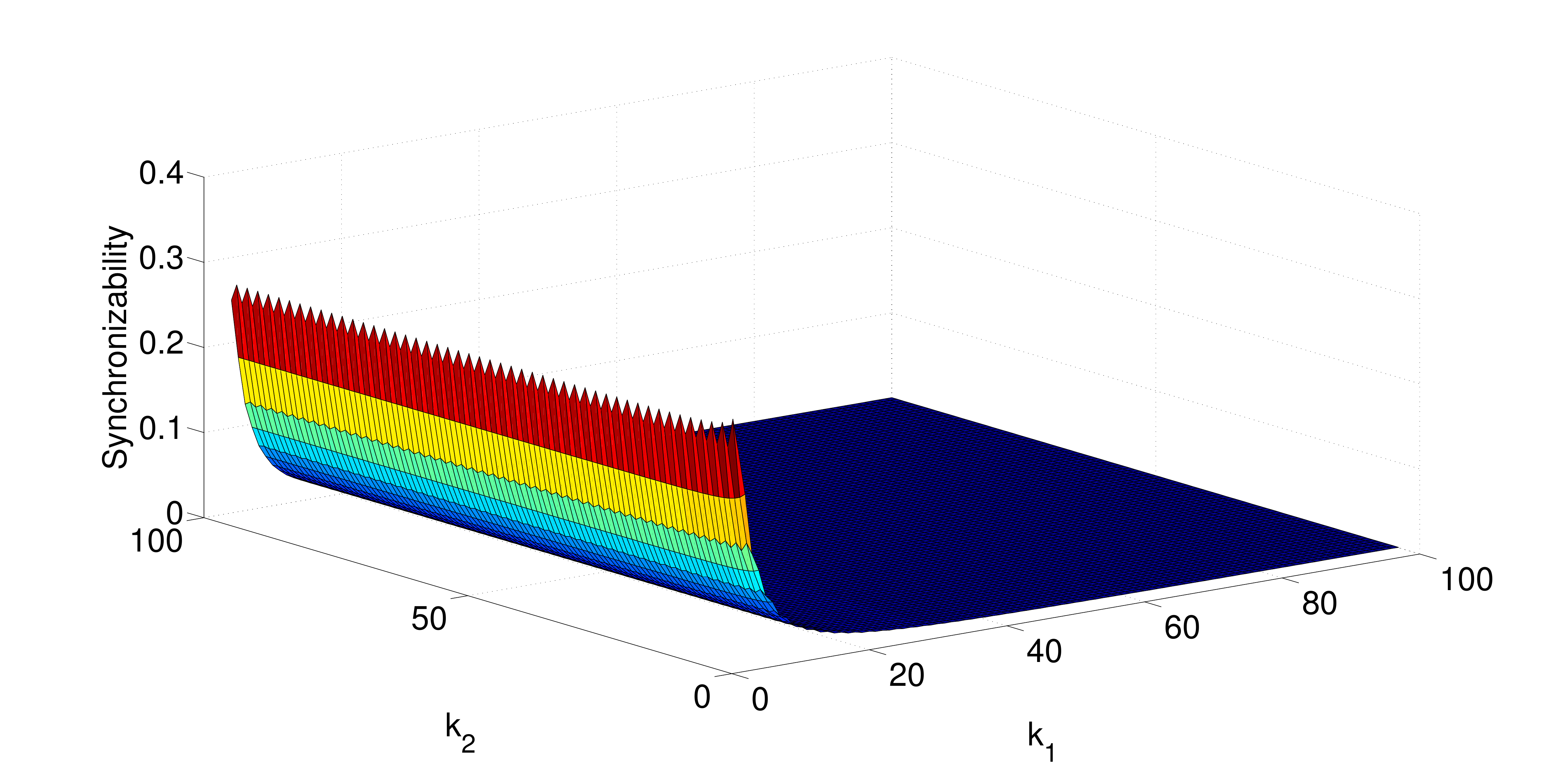}
\caption{Synchronizability versus $k_{1}$ and $k_{2}$ for $r$-nearest neighbor torus}
\label{fig:11}
\end{figure}

\begin{figure}[!t]
\centering
\includegraphics[width=2.5 in]{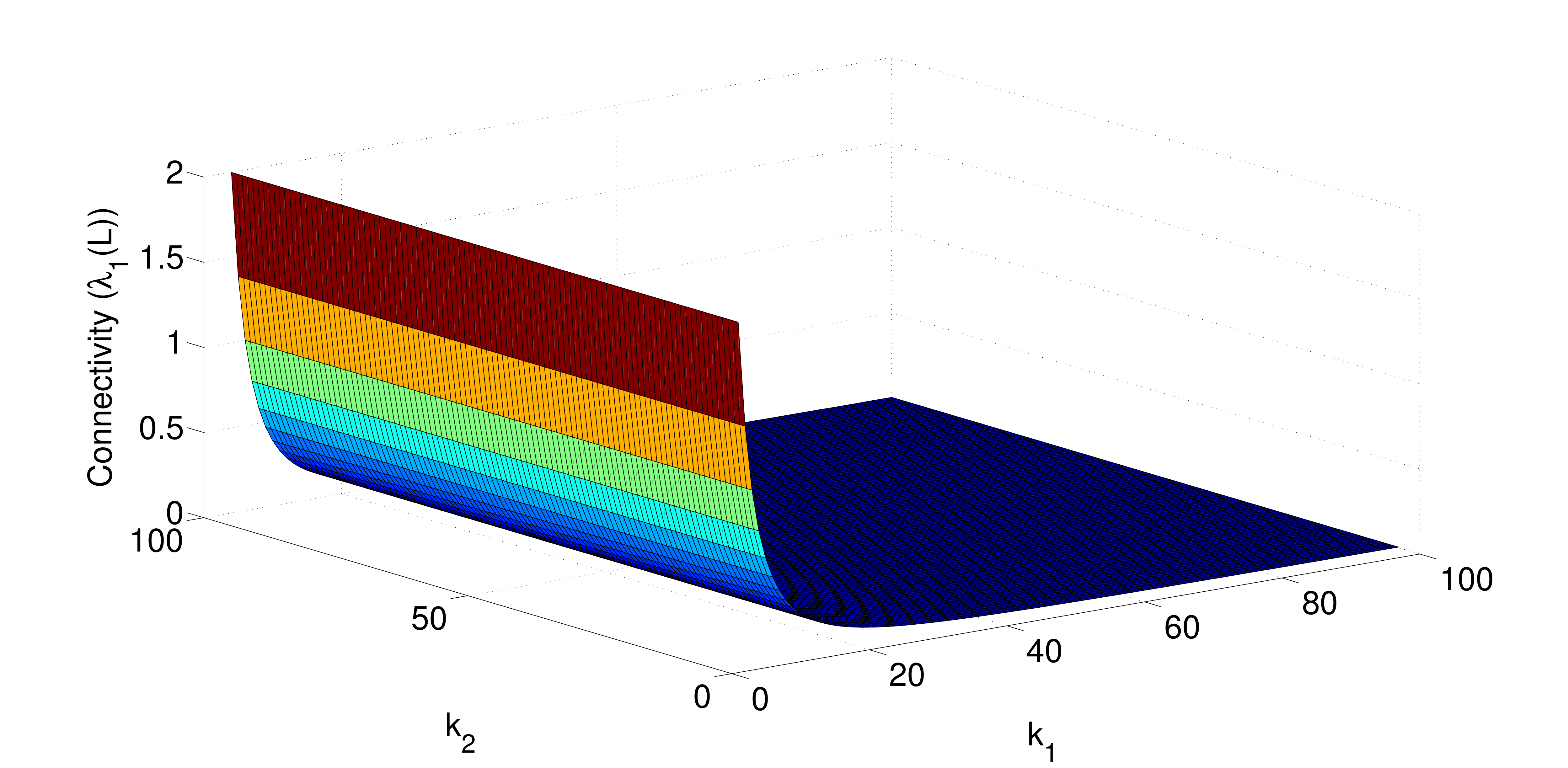}
\caption{Connectivity versus $k_{1}$ and $k_{2}$ for $r$-nearest neighbor torus}
\label{fig:12}
\end{figure}

\begin{figure}[!t]
\centering
\includegraphics[width=2.5 in]{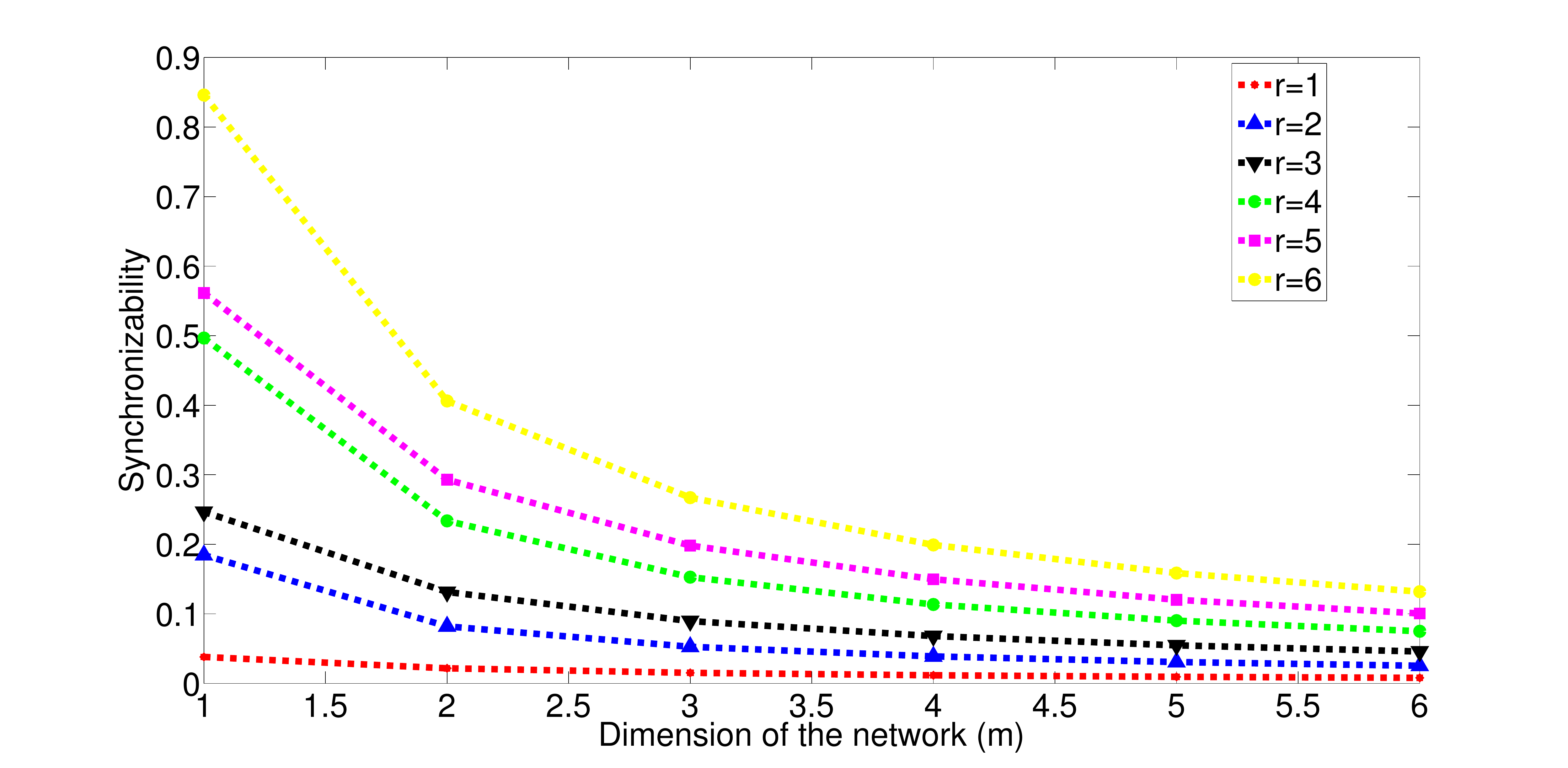}
\caption{Synchronizability versus Dimension for $r$-nearest neighbor torus}
\label{fig:13}
\end{figure}

Network synchronizability $R$ and connectivity $\lambda_{1}(L)$ versus number of nodes $n$ and overhead or nearest neighbors $r$ has been plotted for $r$-nearest neighbor cycle, $r$-nearest neighbor two dimensional and $m$-dimensional torus networks. We compared the analytical expressions derived in the Section IV with simulation results. From Fig. 4 and Fig. 5, we can observe that connectivity $\lambda_{1}(L)$  and synchronizability $R$ increases with overhead $r$. Fig. 6 and Fig. 7 show that connectivity 
$\lambda_{1}(L)$  and synchronizability $R$ decreases with $n$. From Fig. 8 and Fig. 9, we can see that connectivity $\lambda_{1}(L)$  and synchronizability $R$
increases with overhead $r$. We can observe the effect of number of nodes $k_{1}$ and $k_{2}$ on $\lambda_{1}(L)$ and $R$ in Fig. 10 and Fig. 11. $\lambda_{1}(L)$ and $R$ decrease with $k_{1}$ and $k_{2}$. Fig. 12 shows the effect of dimension $m$ on synchronizability $R$ for various $r$ values. Synchronizability $R$ decreases with dimension $m$ and synchronizability is very less for higher dimension and lower overhead values. After dimension $m=5$ values, irrespective of increase in dimension values, $R$ approaches constant synchronizability with the increase in overhead values. From the (\ref{43}), we can see that dimension of the network does not effect the algebraic connectivity.

\section{Conclusions}
Derived the exact formulas for network synchronizability and algebraic connectivity for $r$-nearest neighbor cycle, $r$-nearest neighbor two dimensional torus and $r$-nearest neighbor $m$-dimensional torus networks. We studied the effect of number of nodes and overhead or nearest neighbors effect on network synchronizability and algebraic connectivity. Network synchronizability and algebraic connectivity decrease with the number of nodes and increase with the overhead or nearest neighbors. We also studied the effect of network dimension $m$ on network synchronizability and algebraic connectivity. Network Dimension does not effect the connectivity whereas network synchronizability decreases with the dimension values. The generalized expressions for Laplacian spectra derived in this paper can be used to study the other aspects of the network dynamics. Further, our analysis can also be applied to wireless networks, where the variable $r$ in the $r$-nearest neighbor networks is also captures the transmission radius of the node.

\ifCLASSOPTIONcaptionsoff
  \newpage
\fi

\end{document}